\documentclass[11pt]{amsart}
\usepackage{graphicx}
\usepackage{amssymb}
\usepackage{epstopdf}
\usepackage{amsmath}
\usepackage{amscd}
\usepackage{amsthm,amsfonts,amsxtra,amssymb,amscd}
\usepackage[all]{xy}
\usepackage{enumerate}

\newtheorem{definition}{Definition}
\newtheorem{theorem}{Theorem}[section]

\newtheorem{lemma}[theorem]{Lemma}
\newtheorem{proposition}[theorem]{Proposition}
\newtheorem{corollary}[theorem]{Corollary}

\input xy

\xyoption{all}
\title{Nested sets and Jeffrey Kirwan cycles}

\author{C. De Concini}\address{Dip. Mat. Castelnuovo, Univ. di Roma La
Sapienza, Rome, Italy}\email{deconcin@mat.uniroma1.it}
\author{C. Procesi}\address{Dip. Mat. Castelnuovo, Univ. di Roma La
Sapienza, Rome, Italy}\email{procesi@mat.uniroma1.it}
\thanks{The authors are partially supported by the Cofin 40
\%, MIUR}
\begin{document}
\maketitle
\section{Introduction} In this paper we discuss some new notions  in the
theory of hyperplane  arrangements.   The   paper grew out of our plan to
give an  improved and simplified version of some of the results of Szenes
Vergne \cite{SV}.   

We start from a complex vector space $U$ of finite dimension $r$ and  a
finite central  hyperplane  arrangement in $U^*$, given by a finite set 
$ \Delta\subset U$ of linear equations.  From these data one constructs
the ordered set of subspaces, obtained by intersection of the given
hyperplanes, and the open set $\mathcal A_\Delta$ complement of the union
of the hyperplanes of the arrangement.

This paper consists of 3 parts. Part  1 is a recollection of the results
in \cite{dp}. In part  2   we present 3 new results. The first, of
combinatorial nature,   establishes a canonical bijective correspondence
between the set of no broken circuit bases and  maximal nested sets which
satisfy a condition called properness. 

Next we associate to each proper maximal nested set $\mathbb M$ a
geometric cycle $c_{\mathbb M}$ of dimension $r$ in $\mathcal A_\Delta$.
We show that integration of a top degree differential form over this
cycle is done, by a simple algorithm, taking a multiple residue with
respect to a system of local coordinates.    The last result is the proof
that, under the duality given by integration, the basis of cohomology
given by the forms associated to the no broken circuit bases is dual to
the basis of homology determined by the cycles $c_{\mathbb M}$ . 
  Section 3 is dedicated to the application relevant for the computations
of  \cite{SV}, that is to say the Jeffrey Kirwan cycles.

Finally we wish to thank M. Vergne for explaining to us some of the
theory and  for various discussions and suggestions.

 \subsection{Notations} With the notations of the introduction,  let $U$
be a complex vector space of dimension $r$,    $\Delta\subset U$   a
totally ordered finite set of vectors
$\Delta=\{\alpha_1,\ldots ,\alpha_m\}$. These vectors are the linear
equations of a hyperplane arrangement in $U^*$. For simplicity we also
assume that $\Delta$ spans $U$ and any two distinct elements in $\Delta$
are linearly independent.

An example is a (complete) set of positive roots in a root system,
ordered by any total order which  refines the reverse dominance order.

In the $A_{n-1}$ case we could say that $x_i-x_j\geq x_h-x_k$ if
$k-h\geq j-i$ and, if they are equal if $i\leq h$.

 We want to recall briefly the main points of the theory (cf. \cite{te}).
Let $\Omega_i(\mathcal A_\Delta)$ denote the space of algebraic
differential forms of degree $i$ on $\mathcal A_\Delta$. 

 We shall use implicitly the {\it formality},  that is the fact that the $\mathbb
 Z$
subalgebra of differential forms on $\mathcal A_\Delta$ generated by the
linear forms ${1\over 2\pi i}d\log \alpha, \alpha\in\Delta$ is isomorphic (via De Rham theory) to the integral cohomology of   $\mathcal A_\Delta$.

Formality implies in particular that $\Omega_r(\mathcal
A_\Delta)=H^r\oplus d \Omega_{r-1}(\mathcal A_\Delta),$ for top degree
forms.
$H^r\equiv H^r(\mathcal A_\Delta,\mathbb C)$ is the $\mathbb C$ span of the top degree forms
$\omega_\sigma:=d\log\gamma_1\wedge\dots\wedge d\log\gamma_r$ for all
bases
$\sigma:=\{\gamma_1,\dots,\gamma_r\}$ extracted from $\Delta$.  The forms
$\omega_\sigma$ satisfy a set of linear relations generated by the
following ones. Given $r+1$ elements
$\gamma_i\in\Delta,$ spanning $U,$ we have:
$$\sum_{i=1}^{r+1}(-1)^i d\log\gamma_1\wedge\dots
\check{d\log\gamma_i}\dots\wedge d\log\gamma_r=0.$$ Recall that a {\it
non broken circuit } in $\Delta$ (with respect to the given total
ordering) is an ordered linearly independent subsequence
$\{\alpha_{i_1},\ldots ,\alpha_{i_t}\}$ such that, for each $1\leq
\ell\leq t$ there is no $j<i_{\ell}$ such that the the vectors
$\alpha_j,\alpha_{i_{\ell}},\dots ,\alpha_{i_t}$ are linearly dependent.
In other words $\alpha_{i_\ell}$ is the minimum element of
$\Delta\cap\langle \alpha_{i_\ell},\ldots ,\alpha_{i_t}\rangle$. In
\cite{te} it is proved that the elements
$$({1\over 2\pi i})^r\omega_{\sigma}:=({1\over 2\pi i})^rd\log\gamma_1\wedge\ldots\wedge d\log\gamma_r,$$ where
$\sigma=\{\gamma_1,\ldots,  \gamma_r \}$ runs over all  ordered bases of
$V$ which are non broken circuits, give a linear $\mathbb Z-$basis of the integral cohomology of   $\mathcal A_\Delta$.\subsection{Irreducibles}

Let us now recall some notions from \cite{dp}. Given a subset $S\subset
\Delta$ we  shall denote by $U_S$ the space spanned by $S$.
\begin{definition} Given a subset $S\subset \Delta $ set $\overline
S:=U_S\cap \Delta$, the completion of $S$.  $S$ is called complete if
$S=\overline S$. 

A complete subset $S\subset \Delta $   is called reducible if we can find
a partition $S=S_1\dot{\cup} S_2$, called a {\bf decomposition}  such
that 
$U_S=U_{S_1}\oplus U_{S_2}$, irreducible otherwise.

\end{definition} Equivalently we say that the space $U_S$ is reducible. 
Notice that, in the reducible case $S=S_1\dot{\cup} S_2$,  also $S_1$ and
$S_2$ are complete.

 From this definition it is easy to see \cite{dp}:

\begin{lemma} Given complete sets $A\subset S$ and a decomposition
$S=S_1\dot{\cup} S_2$ of $S$ we have that $A=(A\cap S_1)\dot{\cup}(A\cap 
S_2)$ is a decomposition of $A$.  Let  $S\subset \Delta$ be  complete. 
Then there is a sequence (unique up to reordering) $S_1,\ldots ,S_m$ of
irreducible subsets in $S$ such that
\begin{itemize}
\item $S=S_1\cup\cdots\cup S_m$ as disjoint union. \item
$U_S=U_{S_1}\oplus\cdots\oplus U_{S_m}.$
\end{itemize} The  $S_i$'s are called the irreducible components of $S$
and the decomposition  $S=S_1\cup\cdots\cup S_m$, the irreducible
decomposition of $S$.
\end{lemma}

In the example of root systems, a complete set $S$ is irreducible if and
only  if $S\cup-S$ is an irreducible root system.

 We shall denote by $\mathcal I$ the family of all irreducible subsets in
 $\Delta$.
  \subsection{A minimal model}  In
\cite{dp} we have constructed a  minimal smooth variety  $X_\Delta$
containing
$\mathcal A_\Delta$ as an open set with complement a normal crossings
divisor, plus a proper map $\pi:X_\Delta\to U^*$ extending the identity
of  $\mathcal A_\Delta$. The smooth irreducible components of the
boundary are indexed by the {\it irreducible subsets}. To describe the
intersection pattern between these divisors, in   \cite{dp} we developed
the general theory of nested sets.  Maximal nested sets  correspond
 to special points at infinity, intersections of these boundary divisors.
In the paper \cite{SV},  implicitly the authors use the points at
infinity coming from complete flags which correspond, in the philosophy
of  \cite{dp}, to a {\it maximal}   model with normal crossings. It is
thus not a surprise that by passing from a maximal to a minimal model the
combinatorics gets simplified and the constructions become more canonical.

 Let us recall the
main construction of  \cite{dp}.  For each $S\in\mathcal I$ we have a
subspace $S^{\perp}\subset U^*$ where 
 $S^{\perp}=\{a\in U^*\,|\, s(a)=0,\ \forall s\in S\}.$ We have the
projective space $\mathbb P(U^*/S^{\perp})$ of lines in $U^*/S^{\perp}$ a
map $i:\mathcal A_\Delta\to U^*\times_{S\in\mathcal I}\mathbb
P(U^*/S^{\perp})$. 
 Set $X_\Delta$  equal to the closure of the image  $i(\mathcal A_\Delta)$
in this product. In \cite{dp} we have seen that  $X_\Delta$ is a smooth
variety containing a copy of $\mathcal A_\Delta$ and the complement of
$\mathcal A_\Delta$ in $X_\Delta$ is a union of smooth irreducible
divisors $D_S$, having transversal intersection, indexed by the
elements   $S\in \mathcal I$.

 \subsection{Nested sets}  Still in \cite{dp} we showed that a family
$D_{S_i}$ of divisors indexed by irreducibles $S_i$ has non empty
intersection (which is then smooth irreducible) if and only if the family
is {\it nested} according to:
 \begin{definition} A subfamily
$\mathbb M\subset \mathcal I$ is called nested if, given any subfamily
$\{S_1,\ldots ,S_m\}\subset  \mathbb M$ with the property that for no
$i\neq j$, $S_i\subset S_j$, then
$S:=S_1\cup\cdots\cup S_m$ is complete and the $S_i$'s are the
irreducible components of $S$.
\end{definition}

\begin{lemma} 1) Let   $\mathbb M=\{S_1,\ldots ,S_m\}$ be a nested set.
Then $S:=\cup_{i=1}^mS_i$ is complete. The irreducible components of $S$
are the maximal elements of  $\mathbb M$.

2) Any nested set is the set of irreducible components of the elements of
a flag $A_1\supset A_2\supset \dots\supset A_k$, where each $A_i$ is
complete.
\end{lemma}
\proof 1) By definition of nested set, the maximal elements of $\mathbb
M$ decompose their union which is complete.

2) It is clear that, if $A
\subset B$, the irreducible components of $A$ are contained each in an
irreducible component of $B$. From this follows that the irreducible
components of the sets of a flag form a nested set. Conversely let  
$\mathbb M=\{S_1,\ldots ,S_m\}$ be a nested set.  Set 
$A_1=\cup_{i=1}^mS_i$. Next remove from $\mathbb M$ the irreducible
components of $A_1$ (in $\mathbb M$ by part 1)).  We have a new nested
set to which we can apply the same procedure. Working inductively we
construct a flag of which $\mathbb M$ is the decomposition.\qed

One way of using the previous result is the following, given a basis
$\sigma:\{\gamma_1,\dots,\gamma_r\}\subset \Delta$  one can associate to
$\sigma$ a maximal flag $F(\sigma)$  by setting
$A_i(\sigma):=\Delta\cap\langle\gamma_{i+1},\dots,\gamma_r \rangle$. 
Clearly the map form bases to flags is surjective and from flags to
maximal nested sets is also surjective. We thus obtain a surjective map
from bases to maximal nested sets.  We will see that it induces a
bijiective map between no broken circuit bases and proper maximal nested
sets.

\begin{proposition}  1) Let $A_1\supsetneq A_2\dots\supsetneq A_k$, be a
maximal flag of  complete non empty sets. Then $k=r$ and for each $i$,
$A_i$ spans a subspace of codimension $i-1$.

2)  Let $\Delta=S_1\cup\ldots \cup S_t$ be the irreducible decomposition
of $\Delta$.  

\quad i) Then the $S_i$'s are the maximal elements in $\mathcal I$. 

\quad ii) Every maximal nested set contains each of the elements $S_i$,

\quad\quad  $i=1,\ldots ,t$ and is a union of maximal nested sets in the
sets $S_i$.

3)  Let $ \mathbb M$ be a maximal nested set, $A\in
 \mathbb M$ and $B_1,\dots,B_r\in \mathbb M$   maximal among the elements
in $\mathbb M$ properly contained in $A$. 

Then the subspaces
$U_{B_i}$ form a direct sum and
$$ \dim(\oplus_{i=1}^kU_{B_i})+1=\dim U_A.$$

4) A maximal nested set always has $r$ elements.

\end{proposition}
\proof 1) By definition $A_1=\Delta$ spans $U$. If $\alpha\in
A_i-A_{i+1}$ the completion of $A_{i+1}\cup\{\alpha\}$ must be $A_i$ by
the maximality of the flag. On the other hand by definition $\alpha$ is
not in the subspace spanned by $A_{i+1}$ hence we have that $\dim
U_{A_i}=\dim U_{A_{i+1}}+1$ which implies 1). 2) is immediate from the
definitions.  As for 3) by definition  the subspaces
$U_{B_i}$ form a direct sum and since $A$ is irreducible 
$   \oplus_{i=1}^kU_{B_i} \subsetneq U_A.$  Let $\alpha\in A-\cup_{i=1}^k
B_i$ and $B$ be the completion of $\{\alpha\}\cup \cup_{i=1}^k B_i$. We
must have $B=A$ otherwise we can add the irreducible components of $B$ to
$\mathbb M$ which remains nested, contradicting the maximality. Thus $
\dim(\oplus_{i=1}^kU_{B_i})+1=\dim U_A.$  4) follows from 3) and an easy
induction.
\qed 

A maximal nested set  $\mathbb M$ corresponds thus to a set of $r$
divisors in $X_\Delta$ which by \cite{dp}, intersect transversally in a
single point    $P_{\mathbb M}.$   Let us explicit the example of the
positive roots of type $A_{n-1}$. We think of  such a root as  a pair
$(i,j)$ with $1\leq i<j\leq n$. The  irreducible subsets are   indexed by
subsets (which we display as sequences) $(i_1,\ldots ,i_s)$, $s>1$, with
$1\leq i_1<\cdots <i_s\leq n$. To such a sequence corresponds the set $S$
of pairs supported in the sequence. A family of subsets is nested if any
two of them are either disjoint or one is contained in the other.

In this case a maximal nested set $\mathbb M$ has the following property,
if $A\in \mathbb M$ has $k$ elements and $k>2$ we have two possibilities.
Either the maximal elements of  $\mathbb M$ reduce to one subset with
$k-1$ elements or to two disjoint subsets  $A_1,A_2$ with $A=A_1\cup A_2$.

 We define a map $\phi:\mathcal I\to \Delta$ by
 associating to each $S\in \mathcal I$ its minimum  $\phi(S):=\min(a\in
S)$ with respect to
 the given  ordering.

   For example, in the root system case,  with the ordering given before,
 we have that $\phi(S)$ is the highest root in $S$.

We come to the main new definition:
\begin{definition} A maximal nested set $\mathbb M$ is called proper if
the set  $\phi(\mathbb M)\subset \Delta$ is a basis of $V$.
\end{definition} {\bf Example}\quad  In the $A_n$ case with the previous
ordering
$\phi(i_1,\ldots ,i_s)=x_{i_1}-\nobreak x_{i_s}=(i_1,i_s).$

A proper maximal nested set  $\mathbb M$ is thus encoded by a sequence 
of $n-1$ subsets each having at least two elements,  with the property
that, taking the minimum and maximum for each set, these pairs are all
distinct.

It is  easy to see how to inductively define a bijection between
proper maximal nested sets   and  permutations  of $1,\dots,n$ fixing
$n$.  To see this consider  $\mathbb M$ as a sequence $\{S_1,\ldots,
S_{n-1}\}$ of subsets of $\{1,\ldots ,n\}$ with the above properties.  
We can assume that 
 $S_1=(1,2,\ldots n)$ and  have seen that $\mathbb M':=\mathbb M-\{S_1\}$
has either one or two maximal elements. 
  If $S_2$ is the unique maximal element and $1\notin S_2$, by
 induction we get  a  pemutation $p(\mathbb M')$ of $2,\ldots ,n$. We
then  set $p(\mathbb M)$ equal to the permutation which fixes $1$ and is
equal to $p(\mathbb M')$ on  $2,\ldots ,n$. 
  If $S_2$ is the unique maximal element and $n\notin S_2$, 
 we get  a  pemutation $p(\mathbb M')$ of $1,\ldots ,n-1$. We then  set
$p(\mathbb M)$ equal to the permutation which fixes $n$ and is equal to
$\tau p(\mathbb M')\tau$ on $S_2=\{1,\ldots ,n-1\}$,  $\tau$ being  the
permutation which reverses the order in $S_2$.
   If $S_2$ and $S_3$ are the two maximal elements  so that  $\{1,\ldots
,n\}$ is their disjoint union,  and $1\in S_2$,
$n\in S_3$ then by induction we get two permutations $p_2$ and $p_3$ of
$S_2$ and $S_3$ respectively.   We  then set 
 $p(\mathbb M)$ equal to $p_3$ on $S_3$ and equal to 
$\tau p_2\tau$ on $S_2$ ,  $\tau$ being  the permutation which reverses
the order in $S_2$. In particular this shows that there are $(n-1)!$ 
proper maximal nested sets, which can be recursively constructed. This is
the rank of the top cohomology of the complement of the corresponding
hyperplane arrangement. We will see presently that this is a general
phenomenon.

{\bf Remark} Notice that a proper maximal nested set inherits a total
ordering from the total ordering of $\phi(\mathbb M)$, and that this
ordering is clearly a refinement  of the  partial ordering by reverse
inclusion.

Now fix a maximal nested set $\mathbb M$. We clearly have:
\begin{lemma} Given $\alpha\in \Delta$ there exists a unique minimal
irreducible $S\in \mathbb M$ such that $\alpha\in S$.\end{lemma} This
allows us to define a map $p_{\mathbb M}:\Delta\to \mathbb M$ by setting
$p_{\mathbb M}(\alpha):=S$.
\begin{definition} If $\sigma\subset \Delta$ is a basis of $V$,   we say
that
$\sigma$ is adapted to $\mathbb M$ if the restriction of
$p_{\mathbb M}$ to $\sigma$ is a bijection.
\end{definition}

Notice that, if $\mathbb M$ is proper, then the basis
$\phi(\mathbb M)$ is clearly adapted to $\mathbb M$.

\section{A basis for homology}

We have seen in section 1) that, given a basis
$\sigma=\{\gamma_1,\dots,\gamma_r\}$  we can associate to $\sigma$ a
maximal nested set, which we now denote by $\eta(\sigma)$. 
$\eta(\sigma)$ is the decomposition of the flag
$A_i=\Delta\cap\langle\gamma_i,\dots,\gamma_r\rangle$.   Let us denote by
$\mathcal C$ the set of
   non broken circuit  bases of $V$, by  $\mathcal M$ denote the set of
proper maximal nested set.
 \begin{lemma}\label{circu}  
 If a no broken circuit basis $\sigma$ is adapted to a   proper nested
set $\mathbb M=\{S_1,\ldots ,S_r\},$   then $\sigma=\phi(\mathbb M)$. 
 \end{lemma}

\begin{proof}  
   Let $\sigma=\{\alpha_{i_1},\ldots,\alpha_{i_r}\}$. Clearly $i_1=1$,
and $\alpha_{i_1}$ is the minimum element of $\Delta$. Let $A$ be   the
irreducible component  of $\Delta$   containing
$\alpha_1$.  We have that $A\in \mathbb M, \phi(A)=\alpha_1$ so
$A=S_1$. We claim that
$p_{ \mathbb M}(\alpha_1)=A$.  This follows from the fact that $ \mathbb
M$ is proper so $\alpha_1$ cannot be contained in two distinct elements
$A,B$ of  $\mathbb M$,  otherwise $\phi(A)=\phi(B)$. By  lemma 1.2
$\Delta':=S_2\cup \dots\cup S_r$ is complete. $A\not\subset \Delta'$
 otherwise,
still by 1.2 we would have that $A$ is one of the $S_i,\  i\geq 2$. 
Since $\sigma':=\{\alpha_{i_2},\ldots,\alpha_{i_r}\}$  is adapted to $
\mathbb M':=\{S_2, \dots, S_r\}$ we must have that the space
$U_{\Delta'}$ spanned by $\Delta'$,  is $r-1$ dimensional,
$\{\alpha_{i_2},\ldots,\alpha_{i_r}\}$ is a no broken circuit basis for
$U_{\Delta'}$
relative to  $\Delta'$ ordered by the total order induced from that of
$\Delta$ and  adapted to the proper nested  set $ \mathbb M'$.  We can
thus finish by induction.  
\end{proof}

\begin{theorem}\label{circul}  We have that $\eta$ maps $ \mathcal C$ to
$\mathcal M$ and $\phi$  maps $\mathcal M$ to $ \mathcal C$.
Furthermore $\eta$
and $\phi$ are  bijections which are one the inverse of the other.
 \end{theorem}
 \proof  Let  $\sigma=\{\gamma_1,\dots,\gamma_r\}\in \mathcal C.$ By
definition, for each $i$ we have that $\gamma_i$ is the minimum element
in  $A_i=\Delta\cap\langle\gamma_i,\dots,\gamma_r\rangle$. Thus it is
also the minimum element in one of the irreducibles  decomposing $A_i$.
It follows that  $\eta(\sigma)$ is proper and that
$\phi\eta(\sigma)=\sigma.$  

 Conversely let $\mathbb M=\{S_1,\dots,S_r\}\in \mathcal M$,  and let
$\gamma_i=\phi(S_i)$. By definition the $\gamma_i$'s are linearly
independent, $\gamma_i<\gamma_{i+1}$ and     $\mathbb M$ is the
decomposition of the flag  $A_i:=\cup_{j\geq i}S_j$.  We thus have by the
definition of $\phi$ that $\gamma_i$ is the minimum element in $A_i$.
Since $A_i$ is complete we deduce   that 
$\sigma=\{\gamma_1,\dots,\gamma_r\}\in \mathcal C.$   Clearly 
$\eta(\phi(\mathbb M))=\mathbb M$.\qed 
 \begin{corollary}\label{circu}  
 A no broken circuit basis $\sigma$ is adapted to a   unique maximal
proper nested set $\mathbb M$  and $\sigma=\phi(\mathbb M)$. 
 \end{corollary}

Let us now fix a basis $\sigma\subset\Delta$. Write
$\sigma=\{\gamma_1,\ldots ,\gamma_r\}$ and consider the $r$-form
$$\omega_{\sigma}:=d\log\gamma_1\wedge\ldots\wedge d\log\gamma_r.$$ This
is a holomorphic form on the open set $\mathcal A_{\Delta}$ of $U^*$
which is the complement of the arrangement formed by  the hyperplanes
whose equation is in $\Delta$. In particular if
$\mathbb M\in \mathcal M$, we shall set $\omega_{\mathbb
M}:=\omega_{\phi(\mathbb M)}$.

Also if $\mathbb M\in \mathcal M$, we can define a homology class in
$H_r(\mathcal A_{\Delta},\mathbb
Z)$ as follows. 

Identify $U^*$ with $\mathbb A^r$ using the coordinates $\phi(S),\ S\in
\mathbb M$. Consider another complex affine space $\mathbb A^r$ with
coordinates $z_S$, $S\in \mathbb M$. In $\mathbb A^r$ take the small
torus $T$ of equation
$|z_S|=\varepsilon$ for each $S\in \mathbb M$.  Define a map
$$f:\mathbb A^r\to U^*,\ \ {\rm by} \ \ \phi(S):=\prod_{S'\supset
S}z_{S'}.$$ In \cite{dp} we have proved that this map  lifts, in a
neighborhood of 0, to a local system of coordinates of the model
$X_{\Delta}$  .  To be precise for a vector $\alpha\in\Delta$, set
$B=p_{\mathbb M}(\alpha)$. In the coordinates $z_S$, we have that
\begin{equation}\label{esp}\alpha=\sum_{B'\subset
B}a_{B'}\prod_{S\supseteq B'}z_S=\prod_{S\supseteq
B}z_S(a_B+\sum_{B'\subset B}a_{B'}\prod_{B\supsetneq S\supseteq
B'}z_S)\end{equation} with $a_{B'}\in\mathbb R$ and $a_B\neq 0$.   Set 
$f_{M,\alpha}(z_S):= a_B+\sum_{B'\subset B}a_{B'}\prod_{B\supsetneq
S\supseteq B'}z_S $  and $A_{\mathbb M}$ be the complement in the affine
space $\mathbb A^r$ of coordinates $z_S$ of the hypersurfaces of
equations $f_{M,\alpha}(z_S)=0$. The main point
 is  that, 
$A_{\mathbb M}$ is an open set of $X_\Delta$. The point 0 in  $A_{\mathbb
M}$  is the {\it point at infinity} $P_{\mathbb M}$. The open set
$\mathcal A_\Delta$ is contained in    $A_{\mathbb M}$  as the complement
of   the  divisor with normal crossings given by the equations $z_S=0$.
From this one sees immediately  that,
 if $\varepsilon$ is sufficiently small $f$ maps $T$ homeomorphically
into $\mathcal A_{\Delta}$. Let us give to $T$ the obvious orientation
coming from the total ordering of $\mathbb M$, so that
$H_r(T,\mathbb Z)$ is identified with $\mathbb Z$ and set
$c_{\mathbb M}=f_*(1)\in H_r(\mathcal A_{\Delta},\mathbb Z)$.

\begin{proposition}\label{forme} Let $\sigma=\{\gamma_1,\ldots
,\gamma_r\}\subset \Delta$ be a basis of $V$. Let $\mathbb M\in
\mathcal M$. Then

1) If $\sigma$ is not adapted to $\mathbb M$,
$$\int_{c_{\mathbb M}}\omega_{\sigma}=0.$$

2) If $\sigma$ is  adapted to $\mathbb M$,  consider the sequence
$p_{\mathbb M}(\gamma_1),\ldots ,p_{\mathbb M}(\gamma_r)$. This is a
permutation $\pi$ of the totally ordered set $\mathbb M$ and we denote by
$s(\mathbb M, \sigma)$ its sign. Then
$${1\over  (2\pi i)^r}\int_{c_{\mathbb M}}\omega_{\sigma}=s(\mathbb M,
\sigma).$$
\end{proposition}
\begin{proof} Given $\alpha\in\Delta$, from equation (1) we deduce 
that,  in the neighborhood $A_{\mathbb M}$,  the 1-form $d\log\alpha$  equals the sum of the 1-form  $\sum_{S\supseteq B}d\log z_S$  and of
a 1-form
$\psi_B:=  d\log (a_B+\sum_{B'\subset B}a_{B'}\prod_{B\supsetneq
S\supseteq B'}z_S)$  which is exact and holomorphic on the solid torus in
$\mathbb A^r$ defined by $|z_S|\leq \varepsilon$.

When we substitute these expressions in the linear forms $d\log
\gamma_i$  and expand the product $\omega_{\sigma}$ we obtain various
terms. Some terms vanish since we repeat twice a factor $d\log z_S$, 
some terms contain  a factor 
$\psi_B$ hence they are exact.  The only possible contribution which
gives a non exact form is when $\sigma$ is adapted to $\mathbb M$, and then it is given by  the term  $s(\mathbb M,
\sigma)\omega_{\mathbb M}$. From this observation both 1) and 2) easily
follow.
\end{proof}

Given the class  $c_{\mathbb M}$ and an $r-$dimensional differential form
$\psi$ we can compute $\int_{c_{\mathbb M}}\psi$. Denoting by $P_{\mathbb
M}$ the point at infinity corresponding to 0 in the previously
constructed coordinates $z_i:=z_{S_i}$ we shall say:
\begin{definition} The integral  ${1\over  (2\pi i)^r}\int_{c_{\mathbb
M}}\psi$ is called the {\bf residue} of $\psi$ at the point at infinity
$P_{\mathbb M}$. We will also denote it by $res_{\mathbb M}(\psi)$.
\end{definition}

{\bf NOTE}  The algebraic forms,  in a neighborhood of the point 
$P_{\mathbb M}$ and in the coordinates $z_i$ have the form
$\psi=f(z_1,\dots,z_r)dz_1\wedge\dots\wedge dz_r$ with $f(z_1,\dots,z_r)$
a Laurent series which can be explicitely computed.  Then the residue
$res_{\mathbb M}(\psi)$  equals the coefficient of
$(z_1\dots z_r)^{-1}$,  in this series.\smallskip 

 We can summarize this section with the
main Theorem.

\begin{theorem} The set of elements $c_{\mathbb M}$, $\mathbb M\in
\mathcal M$ is the basis of $H_r(\mathcal A, \mathbb
 Z),$ dual, under the residue
pairing,  to the basis given by the forms
$\omega_{\phi(\mathbb M)}$:  the forms  associated to  the no broken
circuit bases relative to the given ordering.
\end{theorem}
\begin{proof}  This is a consequence of   2.2,  2.3, 2.4.\end{proof}

 {\bf
Remarks.}
1) The formulas found give us an explicit formula for the projection
$\pi$ of $\Omega_r(\mathcal A)=H^r\oplus d
\Omega_{r-1}(\mathcal A_\Delta)$ to $H^r$ with kernel $d
\Omega_{r-1}(\mathcal A_\Delta)$. We have:
\begin{equation}\pi(\psi)=\sum_{\mathbb M\in \mathcal M
}res_{\mathbb M}(\psi)\omega_{\mathbb M}. \end{equation} 

2)  Using the projection $\pi$ any linear map on $H^r$, in particular the
Jeffrey Kirwan residue (see below), can be thought of as a linear map on $\Omega_r(\mathcal A)$
vanishing on $ d
\Omega_{r-1}(\mathcal A_\Delta)$. Our geometric description of homology
allows us to describe any such map as integration on a cycle, linear
combination of the cycles $c_{\mathbb M}$.   

 3) There are several
possible applications of these formulas to combinatorics and counting
integer points in polytopes. The reader is referred to
\cite{WV},\cite{BL}.

4) We have treated only top homology but all homology can be described in
a similar way due to the fact that for each $k$ the $k^{th}$ cohomology
decomposes into the contributions relative to the subspaces of
codimension $k$ and the corresponding transversal configuration.
\section{The Jeffrey-Kirwan residue} In this section $V$ is a real
$r-$dimensional vector space and
 $U:=V\otimes_{\mathbb R}\mathbb C$,
$\Delta=\{\alpha_1,\dots,\alpha_n\}\subset V$.  
 We further restrict to the case in which  there exists a linear function
on $V$ which is positive on $\Delta$.

  Now let us assume that we have  fixed once and for all an orientation
of $V^*$ by choosing an ordered basis $\xi=(x_1,\ldots ,x_r)$ of
$V$ and taking the orientation form $d\underline x=dx_1\wedge
dx_2\wedge\cdots \wedge dx_r$,  for example if $\Gamma$  is the lattice
spanned by the vectors in
$\Delta$  we can take an ordered basis of
$\Gamma$.  

This gives a canonical way of identifying the $r-$forms with functions on
$\mathcal A_\Delta$.   The form
$\omega_\sigma =d\log\gamma_1\wedge\dots\wedge d\log\gamma_r$ is
identified with the function 
$\epsilon_\sigma n_\sigma^{-1}\prod_i\gamma_i^{-1}$ where $n_\sigma$ is
the index in $\Gamma$ of the lattice spanned by the elements $\gamma_i$ (
a volume element), and $\epsilon_\sigma$ a sign expressing the
orientation. Let ${\tilde H^r}$ denote this space of functions.  
Finally, when all the elements in
$\Delta$ are on the same side of some hyperplane there is another
interesting way of representing
$H^r$ in which the Jeffrey Kirwan residue appears in a very natural
way.  This is done via the Laplace transform
$\int e^{-(x,y)}f(y)dy$.  Which in our setting has to be understood as a
transform from functions on $V^*$ with a prescribed invariant Lebesgue
measure (the one induced by $\Gamma$) to functions on $V$ (more
intrinsically to $r$ differential forms).

 Precisely consider the cone   $C$ spanned by the vectors in $\Delta$. For
each basis
$\sigma:=\{\gamma_1,\dots,\gamma_r\}$ extracted from $\Delta$ let
$C(\sigma)$ be the positive cone that it generates and $\chi_\sigma$ its
characteristic function. Let finally $K^r$  be the vector space spanned
by the functions $\chi_\sigma$.  From the basic formula
$\int_0^\infty\dots
\int_0^\infty e^{-\sum_{i=1}^rx_iy_i} dy_1\dots dy_r={1\over x_1\dots
x_r}$ and linear coordinate changes,  it  is   easy to verify that, 
 the Laplace transform of $\chi_\sigma$ is
$n_\sigma^{-1}\prod_i\gamma_i^{-1}$, as function on the dual positive
cone, consisting of all
$x$ such that $(x,y)>0,\ \forall y\in C$. Therefore combining with the
isomorphism of ${\tilde H^r}$ with $H^r$ we have a Laplace transform
$L:K^r\to H^r$ with  
 $L(\chi_\sigma)=\epsilon_\sigma\omega_\sigma.$  

  $L$ is a linear isomorphism   in which it is easy to reinterpret
geometrically the linear relations previously described.  Finally the
Jeffrey Kirwan residue is a linear function $\psi\to J\langle c\,|\, \psi
\rangle$ on $H^r$  depending on a   regular vector $c$. It corresponds to
the linear function defined on $K^r$ which just consists in evaluating
the functions $f$ in $c$.  In other words $  J\langle c\,|\, \psi
\rangle=L^{-1}(\psi)(c).$  By the definition of $K^r$ it is clear that
this linear function depends only on the chamber $C$ in which $c$ lies. 
Our final result is the description of  a geometric cycle $\delta(C)$ 
such that $J\langle c\,|\, \psi \rangle={1\over  (2\pi
i)^r}\int_{\delta(C)}\psi.$    Given an ordered basis $\tau$ of $V$, we
set $\nu_{\tau}$ equal to 1 if the ordered basis $\tau$ has   the same
orientation as
$\xi$, -1 otherwise. For  a proper maximal nested set ${\mathbb M}$ we
denote $\nu_{\phi({\mathbb M})}$ by $\nu_{\mathbb M}$.

For each  basis $\tau\subset \Delta$,  set $C(\tau)=\{x\in V\,|
x=\sum_{\alpha\in\tau}a_{\alpha}\alpha, a_{\alpha}>0\}$.  Set for
simplicity, for a proper maximal nested set ${\mathbb M}$,
$C({\mathbb M}):=C(\phi({\mathbb M}))$. 

A this pont we can recall some facts from  \cite{SV}. Assume that in
$V$ we have a lattice $\Gamma$ which we interpret as the character group
of an $r-$dimensional torus $T_r$. Assume that
$\Delta\subset \Gamma$ is a set of characters.

We have the following sequence of ideas. First of all we use the elements
$\alpha_i,\ i=1,\dots,n$ to construct an $n-$dimensional representation
$Z$ of $T$ direct sum of the 1-dimensional representations with character
$\alpha_i^{-1}$.  Call
$R=\mathbb C[x_1,\dots,x_n]$ the ring of polynomial functions on
$Z$. On $Z$ we have an action of the $n-$dimensional torus $D_n$ of
diagonal matrices and $T$ acts via a homomorphism into $D_n$. Hence the
torus $T$ acts on $R$ and $x_i$ has weight $\alpha_i$. If $\gamma\in
\Gamma$ is a character, define $R(\gamma)$ to be the subspace of $R$ of
weight $\gamma$ with respect to $T$.

We shall denote by $\mathcal R$ the set of regular vectors i.e. the set
of vectors which cannot be written as a linear combination of the
elements in a subset $S$ of $\Delta$ of cardinality smaller than $r$.

Consider a regular vector
$\xi\in \Gamma$ and let $R_\xi=\oplus_{k=0}^\infty R(k\xi)$. 

The following facts are well known \cite{BV}. $R_\xi$ is a finitely
generated
  subalgebra stable under the torus $D_n$. So, if we  grade $R_\xi$ so
that  $R(k\xi)$ has degree $k$ we can consider the projective variety 
$T_\xi:=Proj(R_\xi)$ with a line bundle $\mathcal L$ such that
$H^0(T_\xi, k\mathcal L)=R(k\xi)$. $T_\xi$  is an embedding of the  $n-r$
dimensional torus  $D_n/T$ and the regularity of $\xi$ implies that
$T_\xi$  is an orbifold. The elements $\alpha_i$ index the boundary
divisors of this torus embedding. Thus to each $\alpha_i$ we can
associate a degree 2 cohomology class,  the Chern class of the
corresponding divisor, which is still expressed with the same symbol
$\alpha_i$.  According to the theory developed by Jeffrey Kirwan,
discussed in \cite{BV} one can compute the intersection numbers
$\int_{T_\xi} P(\alpha_1,\dots,\alpha_n)$ using the notion of Jeffrey
Kirwan residue.     Denoting by $C$ the chamber in which $\xi$ lies, one
has:
\begin{equation}\label{SV}\int_{T_\xi}
P(\alpha_1,\dots,\alpha_n)=J\langle c\,|\,{P(\alpha_1,\dots,\alpha_n)\over
\alpha_1\alpha_2\dots\alpha_n}d\mu \rangle .\end{equation}  We want to
represent this residue by integration over a cycle:
\begin{equation}\label{SW} J\langle c\,|\,\psi \rangle= {1\over  (2\pi
i)^r}\int_{\delta(C)}\psi.\end{equation}
  From the formula  $L(\chi_\sigma)=\nobreak\epsilon_\sigma\omega_\sigma
,$ and defininiton of the Jeffrey Kirwan cycle   discussed in the
introduction, we see that  $\delta(C)$  is the $r$-cycle whose value on
the $r$-form $\omega_{\sigma}$, for every basis $\sigma=\{\gamma_1,\ldots
,\gamma_r\}\subset \Delta$  of $V$,  is given by 
\begin{equation}{1\over  (2\pi
i)}\int_{\delta(C)}\omega_{\sigma}=\begin{cases} 0 \text {\ \ if \ }C\cap
C(\sigma)=\emptyset\\ \nu_{\sigma} \text {\ if \ }C\subset C(\sigma)
\end{cases}\end{equation}   Using this description   of $\delta(C)$ and
the fact that our homology  basis $c_{\mathbb M}$ is dual to the
cohomology basis $\omega_{\phi(\mathbb M)}$,  one immediately has:
\begin{theorem}
$$\delta(C)=\sum_{{\mathbb M}\in \mathcal M |C\subset C( {\mathbb
M})}\nu_{\mathbb M} c_{\mathbb M}.$$
\end{theorem}

\end{document}